\input amstex
\documentstyle{amsppt}

\topmatter

\title
A Class of Second Order Tangent Sets
\endtitle
\date   June 16, 2020\enddate
\author      S.~S.~Kutateladze\endauthor
\address     Sobolev Institute of Mathematics\hfill\break
             Novosibirsk 630090, Russia
\endaddress
\email  sskut\@math.nsc.ru\endemail
\thanks
Accepted for publication in 
{\it Siberian Mathematical Journal}\copyright 2020 Pleiades Publ. Ltd.\hfill\break
 http://pleiades.online/
\endthanks


\abstract Under consideration are the construction and properties of some special class
of second other tangent sets on using the technique of nonstandard analysis.
\endabstract
\endtopmatter

\def\dom{\operatorname{dom}}
\def\Cl{\operatorname{Cl}}
\def\cl{\operatorname{cl}}
\def\Ha{\mathop{\fam0 Ha}\nolimits}
\def\Bo{\mathop{\fam0 Bo}\nolimits}
\def\A{\mathop{\fam0 A}\nolimits}
\def\cnt{\mathop{\fam0 cnt}\nolimits}

\document

\noindent
Let  $X$ be a real vector space.
Assume that we are given some almost vector topology
~$\sigma $ with the zero neighborhood filter
$\Cal N _\sigma :=\sigma (0)$ as well as
some almost vector topology
~$\tau $ with the filter
$\Cal N _\tau :=\tau (0)$.

Recall that  every almost vector topology
$\sigma $ on $X$ is characterizes by the
two properties: Firstly,
multiplication by each scalar is continuous;  and,
secondly, addition is jointly continuous.
It is clear that $X$ admits an almost vector topology $\sigma $ such that
$\sigma (0)$ coincides with a fixed filter $\Cal N$
if and only if the monad $\mu(\Cal N)$
is an external vector space over the external field of standard scalars.

In the sequel,
$\sigma $ will be a vector topology, unless stated otherwise explicitly.
It is comfortable to work in the assumption of standard environment within
Nelson internal set theory IST (see~[1]). Recall that the {\it monad} $\mu(\Cal F)$ 
of a standard filter $\Cal F$  is the external intersection of the standard elements of~$\Cal F$.
As usual, introduce the {\it infinite proximity\/}
that is associated with the appropriate uniformity in~ $X$, i.e.
$x_1 \approx {}_\sigma x_2 \leftrightarrow x_1 -x_2 \in\mu (\Cal N _\sigma )$.
Note that the monad   $\mu_\sigma(x)$
of the neighborhood filter  $\sigma (x)$  of the topology $\sigma$
is as follows: $\mu_\sigma(x):= x+\mu(\Cal N_\sigma)$. 
Let $\approx$ stand for the infinite proximity on the reals $\Bbb R$.

Recall that if given are some subset
$F$ of~ $X$ and some point $\bar{x}$ in $X$,
then subdifferential calculus (see~[2] )
deals in particular with the {\it Hadamard, Clarke,} and
{\it Bouligand} cones
$$
\allowdisplaybreaks
\Ha(F,\bar{x}):=
\bigcup\limits_{\Sb U\in\sigma (\bar{x})\\ \lambda >0\endSb}\,
\roman{int}_\tau
\biggl(\,
\bigcap\limits_{\Sb x'\in F\cap U\\ 0<\lambda' \le\lambda \endSb}\,
{F-x'\over\lambda' }\biggr);
$$
$$
\Cl(F,\bar{x}):=
\bigcap\limits_{V^{\vphantom{I^I}}\in\Cal N_\tau  }\,
\bigcup\limits_{\Sb
U\in\sigma (\bar{x})\\
\lambda>0\endSb }\,
\bigcap\limits_{\Sb x'\in F\cap U\\
0<\lambda' \le\lambda\endSb}
\biggl({F-x'\over\lambda' }+V\biggr);
$$
$$
\Bo(F,\bar{x}):=\bigcap\limits_{\Sb
U\in\sigma (x')\\
\lambda>0 \endSb }\ \cl\nolimits_\tau
\biggl(\,
\bigcup\limits_{\Sb x\in F\cap U\\
0<\lambda' \le\lambda\endSb}\,
{F-x'\over \lambda'}\biggr),
$$
where, as usual, $\sigma (\bar{x}):={\bar x}+\Cal N _\sigma  $. If
$h\in\Ha(F,\bar{x})$ then $F$ is often called
{\it epilipshitzian\/} at $\bar{x}$ with respect to~ $h$.
It is obvious that$$
\Ha(F,\bar{x})\subset \Cl(F,\bar{x})\subset \Bo\,(F,\bar{x}).
$$

Considering an extended real function
$f:X\to \Bbb R \cup \{+\infty\}$, the author of~[3] defined the
{\it second order upper subderivative\/} at  $\bar{x}\in X$ along directions $\bar{v}$, $\bar{v}_1$, and  $ \bar{v}_2$ as follows:
 $$
 f_{\bar{x}}^{(2)}(\bar{v}, \bar{v}_1,\bar{v}_2)
 :=\limsup\limits{\Sb (x,\alpha)\rightarrow\bar{x}
 \\\lambda\downarrow 0,\mu\downarrow 0\endSb}
 \inf\limits{\Sb v\rightarrow\bar{v}\\v_{1}\rightarrow\bar{v}_{1}\\
v_{2}\rightarrow\bar{v}_{2}\endSb}O^{2}_{f}
(x,\alpha,\lambda,\mu,v,v_{1},v_{2}),
 $$
 where we use the  limit construction that is attributed to
 Painlev\'e, Kuratowski, and Rockafellar (for instance, see [1, Section 5.3],[4], or~[5])
 while putting
 $$
 \gathered
 O^{2}_{f} (x,\alpha,\lambda,\mu,v,v_{1},v_{2})
\\
:=
\lambda^{-1}\mu^{-1}(f(x+\lambda v_{1}+\mu v_{2}+\lambda\mu v)-
f(x+\lambda v_{1})-f(x+\mu v_{2})+\alpha).
\endgathered
 $$
 Here  $(x,\alpha)\rightarrow\bar{x}$ stands for the convergence to
 $(\bar{x},f(\bar{x}))$
 in the induced topology of the epigraph of~ $f$.
The article [6] contains some approach to the explicit description of the
tangent sets that are determined by  similar constructions. The description
uses the tools of~ IST.

 As we will see soon, it is convenient to slightly modify
 the above construction by inserting the multiplier~4 but
  retaining the previous notations:
  $$
 f_{\bar{x}}^{(2)}(\bar{v}, \bar{v}_1,\bar{v}_2)
 :=\limsup\limits{\Sb (x,\alpha)\rightarrow\bar{x}
 \\\lambda\downarrow 0,\,\mu\downarrow 0\endSb}
 \inf\limits{\Sb v\rightarrow\bar{v}\\v_{1}\rightarrow\bar{v}_{1}\\
v_{2}\rightarrow\bar{v}_{2}\endSb}O^{2}_{f}
(x,\alpha,\lambda,\mu,v,v_{1},v_{2}),
 $$
 where we now put
 $$
 \gathered
 O^{2}_{f} (x,\alpha,\lambda,\mu,v,v_{1},v_{2})
\\
:=
\lambda^{-1}\mu^{-1}(f(x+\lambda v_{1}+\mu v_{2}+4\lambda\mu v)-
f(x+\lambda v_{1})-f(x+\mu v_{2})+\alpha).
\endgathered
 $$

Given  $F\subset X$, denote the {\it indicator function} of~
 $F$ by~ $\delta_{F}$;
 i.~e.,
 $\delta_{F}(x):=0$ at $x\in F$ and $\delta_{F}(x):=\infty$  at~
 $x\notin F$.
 Introduce the set $\Cl^{(2)}(F,\bar{x})(v_{1},v_{2})$ as follows:
$$
v\in\Cl^{(2)}(F,\bar{x})(v_{1},v_{2})\leftrightarrow
(v,v_{1},v_{2})\in\dom(\delta_{F})_{\bar{x}}^{(2)}.
$$

Considering the case of a normed space with $\tau$---the norm topology and
 $\sigma$---the discrete topology, $\Cl^{(2)}(F,\bar{x})(v,v)$ coincides with
 $\A^{(2)} (F,\bar{x},v)$---the {\it  second order attainable direction set\/} to~ $F$
 at~$(\bar{x},v)$ provided that
$\bar{x} + \cnt(v)$ lies in $F$.  As usual, $\cnt(v):=\{\lambda v: \lambda>0, \lambda\approx0\}$
 is the {\it conatus} of~$v$ (see [1, Subsection 5.1.2]). Recall (see, for instance, [7]) that

$$
\A^{(2)} (F,\bar{x},v):=\{h\in X : (\forall   \lambda_n\downarrow 0) (\exists h_n\rightarrow h)\ \bar{x}+\lambda_nv+\frac{1}{2}\lambda^2_nh_n\in F\}.
$$

To simplify bulky formulas we will assume that  $f$ is continuous at $\bar{x}\in F$
with respect to the topology~$\tau$ on~$X$.

\proclaim{Òõþ¨õüð~1}
The following holds:
$$
\gathered
\Cl^{(2)}(F,\bar{x})(v_{1},v_{2})
\\
=\bigcap\limits_{\Sb V\in\Cal{N}_\tau\\
V_{1}\in\sigma(v_{1})\\V_{2}\in\sigma(v_{2})\endSb}
\bigcup\limits_{\Sb x'\in\tau(\bar{x})\\
\lambda_{1}>0\\\lambda_{2}>0\endSb}
\bigcap\limits_{\Sb x'\in F\cap U\\
0<\lambda'\le\lambda_{1}\\0<\lambda''\le\lambda_{2}\endSb}
\bigcup\limits_{\Sb v'\in\frac{F-x'}{\lambda'}\cap V_{1}\\
v''\in\frac{F-x'}{\lambda''}\cap V_{2}\endSb}
\left(\frac{F-x'-\lambda' v'-\lambda'' v''}{4\lambda'\lambda''}+V\right).
\endgathered
$$
\endproclaim

\demo{Proof}
By transfer it suffices to check the case of standard parameters.
Theorem 5.3.11 of [1] yields
$$
\gathered
v\in\Cl^{(2)}(F,\bar{x})(v_{1},v_{2})
\\
\leftrightarrow
(\forall x'\approx_\tau\bar{x},x'\in F)(\forall \lambda'\approx0,
\lambda''\approx0,\lambda'>0,\lambda''>0)
\\
(\exists v'_{1}\approx_\sigma v_{1})(\exists v'_{2}\approx_\sigma v_{2})
(\exists v'\approx_\sigma v)
\\
x'+\lambda' v'_{1}\in F \wedge x'+\lambda'' v'_{2}\in
F\wedge  x'+\lambda' v'_{1}+\lambda'' v'_{2}+4\lambda'\lambda'' v'\in F.
\endgathered
$$
Denote the set on the right-hand side of the claim by~$A$.
Take $v\in\Cl^{(2)}(F,\bar{x})(v_{1},v_{2})$
 and some standard standard neighborhoods $V\in\Cal N_\tau$, $V_{1}\in \sigma(v_{1})$,
and $V_{2}\in\sigma(v_{2})$.
If $\lambda_{1}$ and $\lambda_{2}$ are strictly positive infinitesimal
while  $U$ is an infinitesimal $\tau$-neighborhood of $\bar{x}$, i.~e.
$U\subset\mu_\sigma(\bar{x})$; 
then 
there are some $v'_{1}\approx v_{1}$,
 $v''_{2}\approx_\sigma v_{2}$, and  $v'\approx_\sigma v $  such that
 $x'+\lambda'v'_{1}\in F$, $x'+\lambda''v'_{2}\in F$,
and $x'+\lambda' v'_{1}+\lambda'' v'_{2}+4\lambda'\lambda'' v'\in F$
for all $x'\in F\cap U$, $0<\lambda'\le\lambda_{1}$, and $0<\lambda''\le\lambda_{2}$.
In other words, there exist $v'_{1}\in(F-x')/\lambda'\cap V_{1}$,
$v'_{2}\in(F-x')/\lambda''\cap V_{2}$, and $v'\in v+V$ satisfying the needed
properties. Since the parameters are standard, conclude that
$v\in A$.

Assume now that $v\in A$. Take some standard neighborhoods
 $V\in\Cal N_\tau$, $V_{1}\in\sigma(v_{1})$
 and $V_{2}\in\sigma(v_{2})$ once again.
 By transfer there are  $U\in\sigma(\bar{x})$, $\lambda_{1}>0$,
 and $\lambda_{2}>0$  such that
  $x'+x'v'_{1}+\lambda'' v'_{2}+4\lambda'\lambda'' v'\in F$
  for all  $x'\in F\cap U$, $0<\lambda'\le \lambda_{1}$, and 
  with some    $v'_{1}\in(F-x')/\lambda'\cap V_{1}$,
   $v'_{2}\in(F-x')/\lambda''\cap V_{2}$, and  $v'\in v+V$.
   Recalling that $x'\in U$ if $x'\approx \bar{x}$  and appealing to
   the properties of infinitesimals, we infer by idealization that 
   $v\in\Cl^{(2)}(F,\bar{x})(v_{1},v_{2})$.
   The proof of Theorem 1 is complete.
\enddemo

   \proclaim{Theorem~2}
   The following hold:

{\rm (1)}
If $\Cl^{(2)}(F,\bar{x})(v_{1},v_{2})\ne\varnothing$
then $v_{1}$ and $v_{2}$ belongs to the Clarke cone $\Cl(F,\bar{x})$.

{\rm (2)}
If $v_{1},v_{2}\in\Ha (F,\bar{x})$ then
 $\Cl^{(2)}(F,\bar{x})(v_{1},v_{2})$
 is a closed semigroup.
\endproclaim

 \demo{Proof}
 Claim (1) becomes obvious on recalling that
 in the standard environment we have
 $$
 h\in\Cl(F,\bar{x})\leftrightarrow(\forall x'\approx_\sigma \bar{x},
 x'\in F)(\forall \alpha'>0,\alpha'\approx 0)(\exists h'\approx_\tau h)
 x'+\alpha' h'\in F.
 $$

 Take $u_{1},u_{2}\in\Ha(F,\bar{x})$.
 Without loss of generality, we will proceed in the standard environment.
 Therefore,
 $$
 v_{1}\in\Ha(F,\bar{x})\leftrightarrow
 (\forall x'\approx_\sigma \bar{x},x'\in F)
 (\forall \alpha>0,\alpha\approx 0)(\forall v'_{1}\approx_\tau v_{1})
 x'+\alpha' v'_{1} \in F;
 $$
 $$
 v_{2}\in\Ha(F,\bar{x})\leftrightarrow
 (\forall x'\approx_\sigma \bar{x},x'\in F)
 (\forall \alpha>0,\alpha\approx 0)(\forall v'_{2}\approx_\tau v_{2})
 x'+\alpha' v'_{2} \in F.
 $$
   Assume now that $u_{1},u_{2}\in\Cl^{(2)}(F,\bar{x})(v_{1},v_{2})$.
 By Theorem~1  we can write that
 $$
 \gathered
(\forall x'\approx_\sigma\bar{x},x'\in F)(\forall \lambda'>0,\lambda'\approx 0)
(\forall \lambda''>0,\lambda''\approx 0)
\\
(\exists v'_{1}\approx_\tau v_{1})(\exists v'_{2}\approx_\tau v_{2})
(\exists u'\approx_\tau u_{1})
\\
x'':=x'+\lambda'v'_{1}+\lambda''v'_{2}+4\lambda'\lambda''u'\in F.
\endgathered
 $$
Using the properties of the vector topology $\sigma$ and its monad $\Cal N_\sigma$,
conclude that $x''\approx_\sigma\bar{x}$.
 Recalling Theorem~1 once again, we find  $v''_{1}\approx_\tau v_{1}$;
  $v''_{2}\approx v_{2} $, and $u''\approx u_{2}$
 satisfying $x''\lambda'v''_{1}+\lambda''v''_{2}+4\lambda'\lambda''u''\in F$.
 Put $v':= v'_{1}+v''_{2}$,
 $v'':=v'_{2}+v''_{2}$, and $u:=u'+u''$.
 Undoubtedly, $v'\approx_\tau v_{1}$,
 $v''\approx_\tau v_{2}$, and $u\approx_\tau u_{1}+u_{2}$. Furthermore,
  $x'+\lambda'v'\in F$ and $x'+\lambda''v''\in F$
 since $v_{1}$ and, $v_{2}$ are hypertangents, i.e. elements
 of~$\Ha(F,\bar{x})$. Moreover,
 $$
 \gathered
 x'+\lambda'v'+\lambda''v''+4\lambda'\lambda''u=
 x'+\lambda'v'_{1}+\lambda'v'_{2}+\lambda''v'_{2}
 +\lambda''v''_{2}+4\lambda'\lambda''u'+4\lambda'\lambda''u''
 \\
 =(x'+\lambda'v'_{2}+\lambda''v''_{2}+4\lambda'\lambda''u')+
 \lambda'v'_{2}+\lambda''v''_{2}+4\lambda'\lambda''u''
 \\
 =
 x''+\lambda'v'_{1}+\lambda''v''_{2}+4\lambda'\lambda''u''\in F.
 \endgathered
 $$
 Consequently, $u_{1}+u_{2}\in\Cl^{(2)}(F,\bar{x})(v_{1},v_{2})$.

 To prove closedness, take
  $u_{0}\in\cl_\sigma\Cl^{(2)}(F,\bar{x})(v_{1},v_{2})$
  and some standard neighborhoods $V,V_{1},V_{2}\in\Cal N_\tau$
  such that  $V_{1}+V_{2}\subset V$.
  There is a~standard vector $u\in\Cl^{(2)}(F,\bar{x})(v_{1},v_{2})$
  satisfying $u-u_{0}\in V_{1}$.  Moreover, using Theorem~1, we conclude 
  that there are  $x'\approx_\sigma\bar{x}$, $x'\in F$, $\lambda'\approx0$,
  $\lambda'>0$, $\lambda''\approx 0$, and $\lambda''>0$ such that
   $u'\in u+V_{2}$,  $v'\in v_{1}+W_{1}$, and  $v''\in W_{2}+v_{2}$
 for the previously given standard neighborhoods $W_{1},W_{2}\in\Cal N_\tau$
 satisfying the containments $v'+\lambda'v'\in F$,
   $v'+\lambda''v''\in F$, and
   $x'+\lambda'v'+\lambda''v''+4\lambda'\lambda''u'\in F$.
   This implies easily that $u'\in u+V_{2}\subset u_{0}+V_{1}+V_{2}
   \subset u_{0}+V$.
   By idealization we find
   $v'\approx v_{1}$, $v''\approx_\tau v_{2}$ and $u'_{0}\approx_\tau u_{0}$
   such that $x'+\lambda'v'\in F$, $x'+\lambda''v'\in F$, and
   $x'+\lambda'v'+\lambda''v''+4\lambda'\lambda''u'_{0}\in F$.
   This means that $u_{0}\in\Cl^{(2)}(F,\bar{x})(v_{1},v_{2})$.
\enddemo

\demo{Remark}
   Theorems 1 and 2 can be generalized to the case of the epiderivatives determined from some collection of
   infinitesimals along the lines of [4] and~[8].
\enddemo

\enddocument